\documentclass[10pt,reqno]{amsart}
\usepackage{latexsym}
\usepackage{graphicx}
\usepackage{fancyhdr,amssymb}

\theoremstyle{plain}
\numberwithin{equation}{section}

\theoremstyle{definition}

\begin{document}
\fancyhead{}
\renewcommand{\headrulewidth}{0pt}
\fancyfoot{}
\fancyfoot[LE,RO]{\medskip \thepage}
\fancyfoot[LO]{\medskip The Mathematical Gazette     }
\fancyfoot[RE]{\medskip The Mathematical Gazette}

\setcounter{page}{1}

\title[Another Proof of $e^{x/y}$ being irrational]{Another Proof of $e^{x/y}$ being irrational}
\author{Sourangshu Ghosh}
\address{Department of Civil Engineering\\
                Indian Institute of Technology Kharagpur\\
                Kharagpur,West bengal\\
                India}
\email{sourangshug123@gmail.com   (Corresponding Author)}

\begin{abstract}
    Continued fractions are used to give an alternate proof of $e^{x/y}$ is irrational.
\end{abstract}

\maketitle

\section{Background}
The most well known proof of the irrationality of $e$ was  provided by Joseph Fourier using proof by contradiction \cite{Fourier}. In 1737, Euler gave the first proof of irrationality of $e$ by using the simple continued fraction expansion of $e$ back \cite{Euler1}\cite{Euler2}\cite{Euler3}.

$$e=2+\frac{{1}}{1 + \frac{1}{2 + \frac{1}{1 + \frac{1}{1+\frac{1}{4+\frac{1}{...}}}}}}. \\$$

A simpler proof of this continued fraction was given by Cohn \cite{Cohn}. There are other proofs of $e$ being irrational in the literature such as those by MacDivitt \cite{MacDivitt}, Penesi \cite{Penesi}, Apostol \cite{Apostol}. Higher powers of $e$ were subsequently also proven to be irrational.  The irrationality of $e^{2}$ was proven in \cite{Liouville1} , of  $e^{3}$  in \cite{Hurwitz}, and of $e^{4}$  in \cite{Liouville2}.

A more generalized result where the power is a rational number was proven by Niven in 1985. Ivan Niven \cite{Niven} showed that $e^{x/y}$ is an irrational number using Niven's Polynomials of the form $\frac{x^n(1-x)^n}{n!}$, which can be also be used to prove that $\pi$ is an irrational number. A similar proof was also given by Aigner \cite{Aigner}.

In this document we use continued fraction to give yet another, but simpler, proof of $e^{x/y}$ being irrational, where $x,y$ are integers.

\noindent\rule{0.84in}{0.4pt} \par
\medskip
\indent\indent {\fontsize{8pt}{9pt} \selectfont DOI: 10.35834/YYYY/VVNNPPP \par}
\indent\indent {\fontsize{8pt}{9pt} \selectfont MSC2020: 11J70,11J72 \par}
\indent\indent {\fontsize{8pt}{9pt} \selectfont Key words and phrases: Irrationality, Generalized Continued Fractions,Exponential \par}

\thispagestyle{fancy}

\vfil\eject
\fancyhead{}
\fancyhead[CO]{\hfill Another Proof of $e^{x/y}$ being irrational}
\fancyhead[CE]{S.~Ghosh  \hfill}
\renewcommand{\headrulewidth}{0pt}

\section{Proof}

We will now state the necessary and sufficient condition for the continued fraction proved by Legendre given in Corollary 3, on page 495, in chapter XXXIV on ~~General Continued Fractions" of Chrystal's Algebra Vol.II \cite{Chrystal} to converge into an irrational number.

\noindent \textbf{Theorem:}
\textit{The necessary and sufficient condition that the continued fraction
$$\frac{b_1}{a_1 + \frac{b_2}{a_2 + \frac{b_3}{a_3 + \dots}}}$$}

\textit{ 
\noindent is irrational is that the values $a_{i}, b_{i}$ are all positive integers, and  there is an integer $n$ such that $|a_i| > |b_i|$  for all $i$ greater than $n$.
}

We start with the Continued Fraction Expansion of the hyperbolic tanh function discovered by Gauss \cite{Wall}\cite{Borwein}

$$\tanh z = \frac{z}{1 + \frac{z^2}{3 + \frac{z^2}{5 + \frac{z^2}{...}}}}. \\$$

Then, if we take $z$ to be $\frac{x}{y}$, where $x$ and $y$ are positive integers, we have
$$\tanh z = \frac{\frac{x}{y}}{1 + \frac{(\frac{x}{y})^2}{3 + \frac{(\frac{x}{y})^2}{5 + \frac{(\frac{x}{y})^2}{...}}}}.== \frac{x}{y + \frac{x^2}{3y + \frac{x^2}{5y + \frac{x^2}{...}}}}. \\$$

\noindent \textbf{Corollary:}
\textit{If $x$ and $y$ are integers, then $tan(x/y)$ is irrational
}

Now as $\tanh r=(e^{r}-e^{-r})/(e^{r}+e^{-r})$, it is obvious that if $e^{r}$ is irrational then so is $\tanh r$. Thus $e^{x/y}$ is irrational.

\end{document}